\documentclass[12pt]{article}

\usepackage{theorem}
\usepackage{amssymb}
\usepackage{amsmath}
\usepackage{amsfonts}
\usepackage{times}

\def \RR {\mathbb R}

\def \EE {\mathbb E}

\def \eps {\varepsilon}
\def \vphi {\varphi}

\def \cE {\mathcal E}
\def \cF {\mathcal F}
\def \cN {\mathcal N}
\def \cD {\mathcal D}

\newtheorem{theorem}{Theorem}[section]
\newtheorem{lemma}[theorem]{Lemma}

\newtheorem{proposition}[theorem]{Proposition}
\newtheorem{corollary}[theorem]{Corollary}
 {\theorembodyfont{\rmfamily}}

\begin{document}

\title{Power-law estimates for the central limit theorem
for convex sets}
\author{B. Klartag\thanks{The author is a Clay Research Fellow
and is also supported by NSF grant $\#DMS-0456590$.}}

\date{}
\maketitle

\begin{abstract}
We investigate the rate of convergence in the central limit
theorem for convex sets established in \cite{CLT}. We obtain
bounds with a power-law dependence on the dimension. These bounds
are asymptotically better than the logarithmic estimates which
follow from the original proof of the central limit theorem for
convex sets.
\end{abstract}

\section{Introduction}
\label{intro} This article is a continuation of \cite{CLT}. In
\cite{CLT} we provided a proof for a basic conjecture in convex
geometry (see \cite{ABP}, \cite{BV}), and showed that the uniform
distribution on any high-dimensional convex body has marginals
that are approximately gaussian. Here we improve some quantitative
estimates regarding the degree of that approximation.

\smallskip We denote by $G_{n,\ell}$ the grassmannian of all
$\ell$-dimensional subspaces of $\RR^n$, and let $\sigma_{n,\ell}$
stand for the unique rotationally invariant probability measure on
$G_{n,\ell}$. The standard Euclidean norm on $\RR^n$ is denoted by
$| \cdot |$. For a subspace $E \subseteq \RR^n$ and a point $x \in
\RR^n$ we write $Proj_E(x)$ for the orthogonal projection of $x$
onto $E$. A convex body in $\RR^n$ is a compact, convex set with a
non-empty interior. We write $Prob$ for probability.

\begin{theorem}
Let $1 \leq \ell \leq n$ be integers and let $K \subset \RR^n$ be
a convex body. Let $X$ be a random vector that is distributed
uniformly in $K$, and suppose that $X$ has zero mean and identity
covariance matrix. Assume that $\ell \leq c n^{\kappa}$.

\smallskip Then there exists a subset $\cE \subseteq G_{n,\ell}$ with $\sigma_{n,\ell}(\cE)
\geq 1 - \exp(-c \sqrt{n})$ such that for any $E \in \cE$,
$$ \sup_{A \subseteq E} \left| Prob \{ Proj_E(X) \in A \}
- \frac{1}{(2 \pi)^{\ell/2}} \int_A \exp \left(-\frac{|x|^2}{2}
\right) dx \right| \leq \frac{1}{n^{\kappa}}, $$ where the
supremum runs over all measurable sets $A \subseteq E$. Here, $c,
\kappa > 0$ are universal constants. \label{main_630}
\end{theorem}

Our methods in \cite{CLT} have yielded bounds that depend
logarithmically on the dimension; the estimates in \cite{CLT} are
closer to those in Milman's quantitative theory of Dvoretzky's
theorem (see, e.g., \cite{mil_later} and references therein). In
one of its formulations, Dvoretzky's theorem states that for any
convex body $K \subset \RR^n$ and $\eps > 0$, there exists a
subspace $E \subseteq \RR^n$ of dimension at least $c \eps^2 \log
n$ with
$$ (1 - \eps) \cD \subseteq Proj_E(K) \subseteq (1 + \eps) \cD. $$
Here, $\cD$ is some Euclidean ball in the subspace $E$ that is
centered at the origin, and $c
> 0$ is a universal constant. The logarithmic dependence on the
dimension is known to be tight in Dvoretzky's theorem (consider,
e.g., an $n$-dimensional simplex). In contrast to that, we learn
from Theorem \ref{main_630}  that the uniform measure on $K$, once
projected to subspaces whose dimension is a power of $n$, becomes
approximately gaussian.

\begin{corollary}
Let $1 \leq \ell \leq n$ be integers with $\ell \leq c n^\kappa$.
Let $X$ be a random vector in $\RR^n$ that is distributed
uniformly in some convex body. Then there exist an
$\ell$-dimensional subspace $E \subset \RR^n$ and $r > 0$ such
that
$$ \sup_{A \subseteq E} \left| Prob \{ Proj_E(X) \in A \}
- \frac{1}{(2 \pi r)^{\ell/2}} \int_A \exp \left(-\frac{|x|^2}{2
r} \right) dx \right| \leq \frac{1}{n^{\kappa}}, $$ where the
supremum runs over all measurable sets $A \subseteq E$. Here, $c,
\kappa > 0$ are universal constants.\label{cor_847}
\end{corollary}

Thus, we observe a sharp distinction between the
measure-projection and the geometric-projection of
high-dimensional convex bodies. We did not expect such a
distinction. Our results are also valid for random vectors with a
log-concave density. Recall that a function $f:\RR^n \rightarrow
[0, \infty)$ is called log-concave when
$$ f(\lambda x + (1 - \lambda) y) \geq f(x)^{\lambda}
f(y)^{1-\lambda} $$ for all $x, y \in \RR^n$ and $0 < \lambda <
1$. The characteristic function of any convex set is  log-concave.
As a matter of fact, the assumptions of Theorem \ref{main_630} and
Corollary \ref{cor_847} may be slightly weakened: It is sufficient
for the density of $X$ to be log-concave; the density is not
necessarily required  to be proportional to the characteristic
function of a convex body. A function $f:\RR^n \rightarrow [0,
\infty)$ is isotropic if it is the density of a random vector in
$\RR^n$ with zero mean and identity covariance matrix.

\begin{theorem} Let $n \geq 1$ be an integer and let $X$ be
a random vector in $\RR^n$ with an isotropic, log-concave density.
Then,
$$
 Prob \left \{ \left| \, \frac{|X|}{\sqrt{n}} -1 \, \right| \geq
\frac{1}{n^{\kappa}}
\right \} \leq C \exp \left( -n^{\kappa} \right),
$$ where $C, \kappa > 0$ are universal constants.
  \label{main_408}
\end{theorem}

The bounds we obtain for $\kappa$ from Theorem \ref{main_408} are
not very good. Our proof of Theorem \ref{main_408} works for, say,
$\kappa = 1/15$. See Theorem \ref{prop_408} below for more precise
information. Compare Theorem \ref{main_408} with the sharp
large-deviation estimate of Paouris \cite{pa_cr}, \cite{pa_gafa}:
Paouris showed that under the assumptions of Theorem
\ref{main_408},
\begin{equation}
 Prob \left \{ |X| \geq C  \sqrt{n} \right \} \leq \exp \left(-
 \sqrt{n} \right)
\label{eq_519}
\end{equation}
 for some universal constant $C > 0$. The estimate
(\ref{eq_519}) is known to be essentially the best possible,
unlike the results in this note which probably miss the optimal
exponents.

\medskip   Define
$$ \gamma(t) = \frac{1}{\sqrt{2 \pi}} \exp \left(- \frac{t^2}{2} \right)
\ \ \ \ \ (t \in \RR), $$ the standard gaussian density. We write
$S^{n-1} = \{ x \in \RR^n ; |x| = 1 \}$ for the standard Euclidean
sphere in $\RR^n$. The unique rotationally-invariant probability
measure on $S^{n-1}$ is denoted by $\sigma_{n-1}$. The standard
scalar product in $\RR^n$ is denoted by $\langle \cdot, \cdot
\rangle$. Theorem \ref{main_408} is the basic requirement needed
in order to apply Sodin's moderate deviation estimates
\cite{sodin}. We arrive at the following result:

\begin{theorem}  Let $n \geq 1$ be an integer and let $X$ be
a random vector in $\RR^n$ with an isotropic, log-concave density.

\smallskip Then there exists $\Theta \subseteq S^{n-1}$ with
$\sigma_{n-1}(\Theta) \geq 1 - C \exp(-\sqrt{n})$ such that for
all $\theta \in \Theta$, the real-valued random variable $\langle
X, \theta \rangle$ has a density $f_{\theta}: \RR \rightarrow [0,
\infty)$ with the following two properties:
\begin{enumerate}
\item[(i)] $\displaystyle \int_{-\infty}^{\infty} \left| f_{\theta}(t) - \gamma(t) \right| dt  \leq
\frac{1}{n^{\kappa}}$,
\item[(ii)] For all $|t| \leq n^{\kappa}$ we have
$\displaystyle \left| \frac{f_{\theta}(t)}{\gamma(t)} - 1 \right|
\leq \frac{1}{n^{\kappa}}$.
\end{enumerate}
Here, $C, \kappa > 0$ are universal constants. \label{thm_sodin}
\end{theorem}

The direction of research we pursue in \cite{CLT} and here builds
upon the investigations of Anttila, Ball and Perissinaki
\cite{ABP}, Brehm and Voigt \cite{BV} and others (e.g.,
\cite{bobkov}, \cite{FGP}, \cite{emanuel}). The methods of proof
in this article have much in common with the technique in
\cite{CLT}. As in \cite{CLT}, the basic idea is to show that a
typical multi-dimensional marginal is approximately
spherically-symmetric. Since the marginal is also log-concave,
then most of the mass of the marginal must be concentrated in a
thin spherical shell, and hence the marginal is close to the
uniform distribution on a sphere. An important difference between
the argument here and the one in \cite{CLT}, is the use of
concentration inequalities on the orthogonal group, rather than on
the sphere. Even though the proof here is more technical and
conceptually more complicated than the argument in \cite{CLT}, it
may be considered more ``direct'' in some respects, since we avoid
the use of the Fourier inversion formula.

\medskip As the reader has probably guessed, we write $c, C,
c^{\prime}, \tilde{C}$ etc., and also $\kappa$, for various
positive
 universal constants, whose value may change from one line to the
 next. The symbols $C, C^{\prime},
\bar{C}, \tilde{C}$ etc. denote universal constants that are
assumed to be sufficiently large, while $c, c^{\prime}, \bar{c},
\tilde{c}$ etc. denote sufficiently small universal constants. The
universal constants denoted by $\kappa$ are usually exponents; it
is desirable to obtain reasonable lower bounds on $\kappa$. The
natural logarithm is denoted here by $\log$, and $\EE$ stands for
expectation.

\section{Computations with log-concave functions}

This section contains certain  estimates related to log-concave
functions. Underlying these estimates is the usual paradigm, that
log-concave densities in high-dimension are quite rigid, up to an
affine transformation. We refer the reader, e.g., to \cite[Section
2]{CLT} and to \cite[Section 2]{psitwo} for a quick overview of
log-concave functions and for appropriate references. The
following result of Fradelizi
 \cite[Theorem 4]{fradelizi} will be frequently used.

\begin{lemma} Let $n \geq 1$ and let $f: \RR^n \rightarrow [0, \infty)$ be an integrable,
log-concave function. Denote $x_0 = \int_{\RR^n} x f(x) dx /
\int_{\RR^n} f(x) dx$, the barycenter of $f$. Then, $$ f(x_0) \geq
e^{-n} \sup f. $$ \label{lem_fradelizi}
\end{lemma}

The proof of our next lemma appears in \cite[Corollary 5.3]{CLT}.

\begin{lemma} Let $n \geq 2$ and let $f: \RR^n \rightarrow [0, \infty)$ be an integrable,
log-concave function. Denote $K = \{ x \in \RR^n  ; f(x) \geq
e^{-10 n} \cdot \sup f \}$. Then,
$$ \int_K f(x) dx \geq \left(1 - e^{-n} \right) \int_{\RR^n} f(x) dx. $$
\label{body_K}
\end{lemma}

The following lemma is essentially taken from \cite{psitwo}.
However, the proof in \cite{psitwo} relates only to even
functions. Below we describe a reduction to the even case. Another
proof may be obtained by adapting the arguments from \cite{psitwo}
to the general case. We write $Vol_n$ for the Lebesgue measure on
$\RR^n$, and $D^n = \{ x \in \RR^n ; |x| \leq 1 \}$ for the
centered Euclidean unit ball in $\RR^n$.

\begin{lemma} Let $n \geq 1$ and let $f: \RR^n \rightarrow [0, \infty)$
be an isotropic, log-concave function. Denote $K = \{ x \in \RR^n
; f(x) \geq e^{-n} f(0) \}$. Then, \begin{equation} K \subseteq C
n D^n, \label{eq_1117}
\end{equation}
where $C
> 0$ is a universal constant. \label{lem_300}
\end{lemma}

\emph{Proof:} We may assume that $n \geq 2$ (the case $n=1$
follows, e.g., from \cite[Lemma 3.2]{bobkov}). Suppose first that
$f$ is an even function. Consider the logarithmic Laplace
transform $\Upsilon f(x) = \log \int_{\RR^n} e^{\langle x, y
\rangle} f(y) dy$, defined for $x \in \RR^n$. According to
\cite[Lemma 2.7]{psitwo} the set $T = \{ x \in \RR^n ; \Upsilon
f(x) \leq n \}$ satisfies
\begin{equation}
T \subseteq C n K^{\circ},
\label{eq_224}
\end{equation}
where $K^{\circ} = \{ x \in \RR^n ; \forall y \in K, \langle x, y
\rangle \leq 1 \}$ is the dual body. Since $f$ is isotropic,
$\int_{\RR^n} \langle x, \theta \rangle^2 f(x) dx = 1$ for all
$\theta \in S^{n-1}$. We use Borell's lemma (e.g., \cite[Appendix
III.3]{MS}) and conclude that for any $\theta \in S^{n-1}$,
$$ \int_{\RR^n} \exp \left( c \langle x, \theta \rangle \right) f(x) dx \leq 2. $$
Consequently, $c D^n \subseteq T$. Since $K$ is convex and
centrally-symmetric, the inclusion (\ref{eq_224}) entails that $K
\subseteq C n D^n$. This completes the proof for the case where
$f$ is an even function.

\medskip It remains to deal with the case where $f$ is not necessarily even.
The log-concavity of $f$ implies that $f(\alpha x) / f(0) \leq
(f(x) / f(0))^{\alpha}$ for any $\alpha \geq 1, x \in \RR^n$.
Therefore,
\begin{equation}
 K^{\prime} := \{ x \in \RR^n ; f(x) \geq e^{-10 n} f(0) \} \subseteq 10
 K.
\label{eq_322}
\end{equation}
According to Lemma \ref{lem_fradelizi} and Lemma \ref{body_K},
\begin{equation} e^n f(0) Vol_n(K^{\prime}) \geq \sup f \cdot Vol_n(K^{\prime})
\geq \int_{K^{\prime}} f(x) dx \geq 1/2. \label{eq_1115}
\end{equation} From (\ref{eq_322}) and (\ref{eq_1115}) we see that
\begin{equation}
Vol_n(K) > \frac{c^n}{f(0)}.
\label{eq_326}
\end{equation}
Denote
$$ g(x) = 2^{n/2} \int_{\RR^n} f(y) f(y + \sqrt{2} x) dy. $$
Then $g$ is an even, isotropic, log-concave density on $\RR^n$, as
follows from the Pr\'ekopa-Leindler inequality (for
Pr\'ekopa-Leindler, see, e.g., the first pages of \cite{pisier}).
Moreover,
\begin{equation}
 g(0) \leq 2^{n/2} \sup f \int_{\RR^n}  f(y) dy = 2^{n/2} \sup f \leq (\sqrt{2} e)^n f(0)
\label{eq_334}
\end{equation}
by Lemma \ref{lem_fradelizi}. The set $K$ is convex with $0 \in
K$. For any $x,y \in K / 4$ we have $y + \sqrt{2} x \in K$ and
hence $f(y), f(y + \sqrt{2} x) \geq e^{-n} f(0)$. Therefore, by
(\ref{eq_326}) and (\ref{eq_334}), if $x \in K /4$,
$$ g(x) \geq \int_{K/4} f(y) f(y + \sqrt{2} x) dy
\geq Vol_n(K/4) (e^{-n} f(0))^2  \geq c^n f(0) \geq \tilde{c}^n
g(0). $$ We conclude that
\begin{equation}
 K / 4 \subseteq \{ x \in \RR^n ; g(x) \geq e^{-Cn} g(0) \}
\subseteq C \{ x \in \RR^n ; g(x) \geq e^{-n} g(0) \},
\label{eq_335}
\end{equation}
where the last inclusion follows directly from the log-concavity
of $g$. Recall that $g$ is even, isotropic and log-concave, and
that (\ref{eq_1117}) was already proven for functions that are
even, isotropic and log-concave. We may thus assert  that
\begin{equation}
 \{ x \in \RR^n ; g(x) \geq e^{-n} g(0) \} \subseteq C n D^n.
\label{eq_712}
\end{equation}
The conclusion of the lemma  thereby follows from (\ref{eq_335})
and (\ref{eq_712}). \hfill $\square$

\begin{corollary} Let $n \geq 1$ and let $f :\RR^n \rightarrow [0, \infty)$
be an isotropic, log-concave function. Then, for any $x \in
\RR^n$,
$$ f(x) \leq f(0) \exp \left( C n - c |x| \right) $$
for some universal constants $c, C > 0$.
\label{lem_309}
\end{corollary}

\emph{Proof:} According to Lemma \ref{lem_300}, for any $x \in
\RR^n$,
$$ f(x) \leq e^{-n} f(0)  \ \ \ \text{if} \ \ |x| \geq C n. $$
By log-concavity, whenever $|x| \geq Cn$,
$$ e^{-n} f(0)   \geq f \left( \frac{C n}{|x|} \cdot x\right) \geq f(0)^{1 - C n /|x|}
\cdot f(x)^{Cn / |x|}. $$ Equivalently,
\begin{equation}
 f(x) \leq f(0) e^{-|x| /C } \ \ \ \text{whenever} \ \ |x| \geq C n.
\label{eq_202}
\end{equation}
According to Lemma \ref{lem_fradelizi} we know that $f(x) \leq e^n
f(0)$ for any $x \in \RR^n$. In particular, $f(x) \leq f(0) \exp(2
n - |x| / C)$ when $|x| < C n$. Together with (\ref{eq_202}) we
obtain
$$
 f(x) \leq f(0) e^{2n - |x| / C} \ \ \ \text{for all} \ \  x \in
\RR^n. $$ This completes the proof. \hfill $\square$

\medskip
We will make use of the following elementary result.
\begin{lemma} Let $f: [0, \infty) \rightarrow [0, \infty)$ be
a measurable function with $0 < \int_0^{\infty} (1 + t) f(t) dt <
\infty$. Suppose that $g: [0, \infty) \rightarrow (0, \infty)$ is
a monotone decreasing function. Then,
$$ \frac{\int_0^{\infty} t^2 f(t) g(t) dt}{\int_0^{\infty} f(t) g(t)
dt} \leq \frac{\int_0^{\infty} t^2 f(t) dt}{\int_0^{\infty} f(t)
dt}.
$$
\label{lem_738}
\end{lemma}

\emph{Proof:} Since $g$ is non-increasing, for any $t \geq 0$,
$$ \frac{\int_0^t f(s) g(s) ds}
{\int_t^{\infty} f(s) g(s) ds} \geq \frac{g(t) \int_0^t f(s) ds}
{g(t) \int_t^{\infty} f(s) ds} = \frac{\int_0^t f(s)
ds}{\int_t^{\infty} f(s) ds}. $$ Equivalently, for any $t \geq 0$,
$$ \frac{\int_t^{\infty} f(s) g(s) ds}
{\int_0^{\infty} f(s) g(s) ds} \leq \frac{\int_t^{\infty} f(s) ds}
{\int_0^{\infty} f(s)  ds}. $$ We conclude that
$$ \frac{\int_0^{\infty} t^2 f(t) g(t) dt}{\int_0^{\infty} f(t) g(t)
dt} = \frac{\int_0^{\infty} \int_{\sqrt{t}}^{\infty} f(s) g(s) ds
dt} {\int_0^{\infty} f(t) g(t) dt} \leq \frac{\int_0^{\infty}
\int_{\sqrt{t}}^{\infty} f(s) ds dt} {\int_0^{\infty} f(t)  dt} =
\frac{\int_0^{\infty} t^2 f(t) dt}{\int_0^{\infty} f(t) dt}.
$$
\hfill $\square$

\medskip  The
identity matrix is denoted here by $Id$. For a $k$-dimensional
subspace $E \subseteq \RR^n$ and for $v
> 0$ we define $\gamma_{E}[v]: E \rightarrow [0, \infty)$ to be
the density $$ \gamma_{E}[v](x) = \frac{1}{(2 \pi v)^{k/2}} \exp
\left( -\frac{|x|^2}{2v} \right) \ \ \ \ \ \ \ \ \ \ \ (x \in E).
$$ That is, $\gamma_{E}[v]$ is the density of a gaussian random
vector in $E$ with mean zero and covariance  matrix that equals $v
Id$. A standard gaussian random vector in $E$ is a random vector
whose density is $\gamma_{E}[1]$. We abbreviate $\gamma_n[v]$ for
$\gamma_{\RR^n}[v]$.

\begin{corollary}  Let $n \geq 1$ and let $f: \RR^n \rightarrow [0, \infty)$ be an
isotropic function. Let $x \in \RR^n, v > 0$ and denote
$$ g(y) = f(y) \exp \left( -\frac{|x-y|^2}{2 v} \right) \ \ \ \ \ \ \ \ (y \in \RR^n), $$
an integrable function. Suppose that $X$ is a random vector in
$\RR^n$ whose density is proportional to $g$. Then,
\begin{enumerate}
\item[(i)] $\displaystyle \EE |X - x|^2 \leq  n + |x|^2$,
\item[(ii)] $\displaystyle |\EE X - x| \leq \sqrt{n} + |x|
$.
\end{enumerate}
\label{lem_802}
\end{corollary}

\emph{Proof:} For $r > 0$ denote
$$ h(r) = r^{n-1} \int_{S^{n-1}} f(x + r \theta) d\theta $$
where the integration is with respect to the surface area measure
on $S^{n-1}$. By integration in polar coordinates,
$$ \int_0^{\infty} h(r) dr = \int_{\RR^n} f = 1. $$ Let $Y$ be a
random vector in $\RR^n$ that is distributed according to the
density $f$. We integrate in polar coordinates and obtain
\begin{eqnarray} \label{eq_1008}
\lefteqn{ n + |x|^2 = \EE |Y - x|^2 } \\ & = & \int_{\RR^n} |y|^2
f(y+x) dy = \int_0^{\infty} \int_{S^{n-1}} r^{n+1} f(x + r \theta)
d \theta dr = \int_{0}^{\infty} r^2 h(r) dr. \nonumber
\end{eqnarray}
Similarly,
\begin{equation}
 \EE |X - x|^2 = \frac{\int_0^{\infty} r^2 h(r) \exp(-r^2 / 2v)
dr}{\int_0^{\infty} h(r) \exp(-r^2 / 2v) dr}. \label{eq_1009}
\end{equation}
We apply the elementary Lemma \ref{lem_738}, based on
(\ref{eq_1008}) and (\ref{eq_1009}). This proves (i). The
inequality (ii) follows from (i) by Jensen's inequality. \hfill
$\square$

\medskip
The crude estimates in the next two lemmas are the main results of
this section. Our first lemma does not use the log-concavity
assumption in an essential way.

\begin{lemma} Let $n \geq 1$, let $f: \RR^n \rightarrow [0, \infty)$ be an
isotropic, log-concave function, and let $ v > 0$. Denote $g = f
* \gamma_n[v]$, the convolution of $f$ with $\gamma_n[v]$.
Then, for any $x \in \RR^n$ with $|x| \leq 10 \sqrt{n}$,
$$ \left| \nabla \log g(x) \right| \leq C \sqrt{n} / v, $$ where $C > 0$ is a
universal constant. \label{lem_349}
\end{lemma}

\emph{Proof:} Since
$$ g(x) = (2 \pi v)^{-(n/2)}
\int_{\RR^n} f(y) \exp \left( -|x - y|^2 / (2v) \right) dy, $$ we
may differentiate under the integral sign and obtain that
$$ \nabla g(x) = (2 \pi v)^{-(n/2)}
\int_{\RR^n} f(y) \frac{y-x}{ v} \exp \left( -|x - y|^2 / (2v)
\right) dy.  $$ Fix $x \in \RR^n$. Denote $g_x(y) = f(y) \exp
\left( -|y-x|^2 / (2v) \right)$. Then,
\begin{equation} \nabla \log g(x) = \frac{\nabla g(x)}{g(x)} =
v^{-1} \left. \int_{\RR^n} (y - x) g_x(y) dy \right/ \int_{\RR^n}
g_x(y) dy. \label{eq_633}
\end{equation}
Let $X$ be a random vector in $\RR^n$ whose density is
proportional to $g_x$. Then $|\EE X - x| \leq \sqrt{n}  + |x|$, as
we learn from Corollary \ref{lem_802}(ii). We conclude from
(\ref{eq_633}) that
$$ |\nabla \log g(x)| = v^{-1} | \EE X - x | \leq \frac{
\sqrt{n} + |x|}{v} \leq 11 \frac{\sqrt{n}}{v}, $$ since $|x| \leq
10 \sqrt{n}$. The lemma is thus proven. \hfill $\square$

\medskip For a subspace $E \subseteq \RR^n$ we write $E^{\perp} = \{
x \in \RR^n ; \forall y \in E, \langle x, y \rangle = 0 \}$ for
its orthogonal complement. For an integrable function $f: \RR^n
\rightarrow [0, \infty)$, a subspace $E \subseteq \RR^n$ and a
point $x \in E$ we write
$$ \pi_E(f)(x) = \int_{x + E^{\perp}} f(y) dy. $$
That is, $\pi_E(f): E \rightarrow [0, \infty)$ is the marginal of
$f$ onto $E$.

\begin{lemma} Let $n \geq 1, 0 < v \leq 1, x \in \RR^n$ and $e \in
S^{n-1}$.
 Let $f: \RR^n \rightarrow [0, \infty)$ be an isotropic,
log-concave function and denote $g = f * \gamma_{n}[v]$. For
$\theta \in \RR^n$ with $|\theta - e| < 1/2$ define
$$ F(\theta) = \log \int_{-\infty}^{\infty} g(x + t \theta) dt. $$
Assume that $|x| \leq 10 \sqrt{n}$. Then,
$$ \left| \nabla F(e) \right|  \leq C n^{3/2} / v^2, $$ where $C > 0$ is a universal
constant. \label{rotate}
\end{lemma}

\emph{Proof:} For $\theta \in \RR^n$ with $|\theta - e| < 1/2$, denote
$$
G(\theta) = \int_{-\infty}^{\infty} g(x + t \theta) dt = (2 \pi
v)^{-n/2} \int_{\RR^n} \int_{-\infty}^{\infty} f(y) \exp \left( -
\frac{|x + t \theta - y|^2}{2 v} \right) dt dy $$ $$ =
 (2 \pi v)^{-(n-1)/2} \cdot |\theta|^{-1} \cdot \int_{\RR^n} f(y) \exp \left(-
\frac{|Proj_{\theta^{\perp}}(x - y)|^2}{2 v} \right) dy
$$
where $\theta^{\perp} =  \{ x \in \RR^n ; \langle x, \theta
\rangle = 0 \}$. We may differentiate under the integral sign
(recall that $f$ decays exponentially fast at infinity by, e.g.,
\cite[Lemma 2.1]{psitwo}) and obtain
\begin{eqnarray*}
\nabla G(e) = \frac{1}{v (2 \pi v)^{n/2} } \int_{\RR^n} f(y)
\int_{-\infty}^{\infty}
 (y - x - t e) t  \exp \left( -\frac{|x + t e - y|^2}{2v} \right)
dt dy.
\end{eqnarray*}
We write $y - x = z + r e$ where $r \in \RR$ and $z \in
e^{\perp}$. A direct computation reveals that
$$
\int_{-\infty}^{\infty}
 [z + (r - t) e] t  \exp \left( -\frac{|(t - r) e - z|^2}{2v} \right)
dt = \sqrt{2 \pi v} \exp \left( -\frac{|z|^2}{2 v} \right) \cdot
(r z - v e).
$$
Denote $H = e^{\perp}$ and set $g_x(y) =  f(y) \exp
\left(-|Proj_H(x-y)|^2 / (2 v) \right)$ for $y \in \RR^n$. Then,
$$ |\nabla G(e)| \leq
\frac{1}{v (2 \pi v)^{(n-1)/2}} \int_{\RR^n} \left( |Proj_H(y-x)|
\cdot |\langle y-x, e \rangle| + v\right) g_x(y) dy. $$ According
to the Cauchy-Schwartz inequality,
\begin{eqnarray}
 \label{eq_634}
 \lefteqn{|\nabla F(e)| = \frac{|\nabla G(e)|}{G(e)}
 \leq 1 + \frac{\int_{\RR^n}
|Proj_H(x-y)| \cdot |\langle x-y, e \rangle| \cdot g_x(y) dy}{v
\int_{\RR^n} g_x(y) dy} }
\\ & \leq & 1 + \frac{1}{v}  \left( \frac{\int_{\RR^n}
|Proj_H(x-y)|^2 g_x(y) dy}{\int_{\RR^n} g_x(y) dy} \cdot
\frac{\int_{\RR^n} \langle x-y, e \rangle^2 g_x(y)
dy}{\int_{\RR^n} g_x(y) dy} \right)^{1/2}. \nonumber
\end{eqnarray}
Our derivation of the inequality (\ref{eq_634}) relies only on
integrability properties of $f$. The log-concavity and
isotropicity of $f$ will come into play next.
 Let $Y$ be a random vector in $\RR^n$ whose density is
proportional to $g_x$. Then the density of the random vector
$Proj_H(Y)$ is proportional to
$$ y \mapsto \pi_H(f)(y) \cdot \exp \left(-\frac{|y - Proj_H(x)|^2}{2 v}
\right) \ \ \ \ \ (y \in H). $$ The density $\pi_H(f)$ is
isotropic and log-concave, by Pr\'ekopa-Leindler. We may thus
apply Corollary \ref{lem_802}(i) and conclude that
$$
 \EE |Proj_H(Y) - Proj_H(x)|^2 \leq (n-1) + |Proj_H(x)|^2. $$
Therefore,
\begin{equation}
\frac{\int_{\RR^n} |Proj_H(y-x)|^2 g_x(y) dy}{\int_{\RR^n} g_x(y)
dy} \leq 2 \max \left \{ n, |x|^2 \right \}. \label{eq_847}
\end{equation}
Next, we deal with the second factor of the product in
(\ref{eq_634}). We write $H^+ = \{ y \in \RR^n ; \langle y, e
\rangle \geq 0 \}$ and $H^- = \{ y \in \RR^n ; \langle y, e
\rangle \leq 0 \}$. Let $g^+: \RR^n \rightarrow [0, \infty)$ be
the function defined by
$$ g^+(y) = \left( \langle y - x, e \rangle + 2|x| + 1 \right)^2 g_x(y) \ \ \ \text{for} \ y \in
H^+ $$ and $g^+(y) = 0$ for $y \not \in H^+$.
 Observe that $g^+$ is log-concave, since both
 $y \mapsto (\langle y - x, e \rangle + 2|x| + 1)^2$ and $y \mapsto
g_x(y)$ are log-concave on $H^+$. Additionally, $g^+$ is
integrable, since $f$ decays exponentially fast at infinity. We
claim that
\begin{equation}
g^+(y) < e^{-10 n} g^+(0) \ \ \ \text{if} \ |y| >  \tilde{C}
\left( \frac{|x|^2}{v} + |x| + n \right). \label{eq_905}
\end{equation}
Indeed, by Corollary \ref{lem_309},  $f(y) \leq f(0) e^{C n - c
|y| }$ for all $y \in \RR^n$. Hence, when $|y| > \tilde{C}
 \left( |x|^2 / v + |x| + n \right)$,
 \begin{eqnarray*}
 \lefteqn {g^+(y) \leq (1 + 2|x| + |y-x|)^2 g_x(y) \leq (1 + 2 |y|)^2
 \cdot f(0) e^{C n - c |y| } } \\ & \leq &
f(0) e^{C^{\prime} n - c^{\prime} |y|}  <  f(0) \exp \left(-10 n -
\frac{|x|^2}{2 v} \right) \leq e^{-10 n} g_x(0) \leq e^{-10 n}
g^+(0),
\end{eqnarray*}
and (\ref{eq_905}) is proven. Denote $K^+ = \{ y \in \RR^n ;
g^+(y) \geq e^{-10 n} g^+(0) \}$. Then, by (\ref{eq_905}) and by
Lemma \ref{body_K},
\begin{eqnarray*}
\lefteqn{ \left[ 1 + 3|x| + \tilde{C} \left( \frac{|x|^2}{v} + |x|
+ n \right) \right]^2 \cdot \int_{K^+} g_x(y) dy } \\ & \geq &
\int_{K^+} g^+(y) dy \geq \left(1 - e^{-n} \right) \int_{\RR^n}
g^+(y) dy. \end{eqnarray*} We deduce that
\begin{equation}
\frac{\int_{H^+} \langle y - x, e \rangle^2 g_x(y)
dy}{\int_{\RR^n} g_x(y) dy} \leq \frac{\int_{\RR^n} g^+(y)
dy}{\int_{K^+} g_x(y) dy} \leq C \max \left \{ \frac{|x|^4}{v^2},
|x|^2, n^2 \right \}. \label{eq_1103}
\end{equation}
The proof that \begin{equation} \int_{H^-} \langle y - x, e
\rangle^2 g_x(y) dy \leq C \max \left \{ \frac{|x|^4}{v^2}, |x|^2,
n^2 \right \} \cdot \int_{\RR^n} g_x(y) dy \label{eq_1104}
\end{equation}
is completely analogous. One just needs to work with the function
$y \mapsto g^-(y) = \left( \langle y - x, e \rangle - 2|x| -
1\right)^2 g_x(y)$, which is log-concave on $H^-$. By adding
(\ref{eq_1103}) and (\ref{eq_1104}) we find that
\begin{equation}
\frac{\int_{\RR^n} \langle y - x, e \rangle^2 g_x(y)
dy}{\int_{\RR^n} g_x(y) dy} \leq C^{\prime} \max \left \{
\frac{|x|^4 + 1}{v^2}, n^2 \right \}, \label{eq_1105}
\end{equation}
since $0 < v \leq 1$. We combine  (\ref{eq_1105}) with
(\ref{eq_847}) and (\ref{eq_634}), and conclude that
$$ |\nabla F(e)| \leq 1 + \frac{C}{v^2} \max \{ \sqrt{n}, |x| \} \cdot
\max \{ |x|^2 + 1, n v \}. $$ The lemma follows, since $|x| \leq
10 \sqrt{n}$ and $0 < v \leq 1$. \hfill $\square$

\section{Concentration of measure on the orthogonal group}

Let $f: \RR^n  \rightarrow [0, \infty)$ be an integrable function,
let $E \subseteq \RR^n$ be a subspace and let $x \in E$. Recall
that we define
$$ \pi_E(f)(x) = \int_{x + E^{\perp}} f(y) dy. $$
We will consider the group $SO(n)$, that consists of all
orthogonal transformations in $\RR^n$ of determinant one. The
group $SO(n)$ admits a canonical Riemannian metric which it
inherits from the obvious embedding $SO(n) \subset \RR^{n^2}$
(that is, a real $n\times n$-matrix has $n^2$ real numbers in it).
 For $U \in SO(n)$ we set
\begin{equation}
 M_{f, E, x}(U) = \log \pi_{E}(f \circ U)(x) = \log \pi_{U(E)}(f) (U x). \label{eq_1141}
\end{equation}
Clearly, for any $U_1, U_2 \in SO(n)$,
\begin{equation}
 M_{f, E, x}(U_1 U_2) = M_{f,U_2(E),U_2(x)}(U_1).
\label{eq_642}
\end{equation} For $U_1, U_2 \in SO(n)$ we write $d(U_1, U_2)$ for the
geodesic distance between $U_1$ and $U_2$ in the connected
Riemannian manifold $SO(n)$.  It is well-known that for any $U_1,
U_2 \in SO(n)$,
$$ \| U_1 - U_2 \|_{HS} \leq d(U_1, U_2) \leq \frac{\pi}{2} \| U_1 - U_2
\|_{HS} $$ where $\| \cdot \|_{HS}$ stands for the Hilbert-Schmidt
norm, i.e., for a matrix $A = (A_{i,j})_{i,j=1,...,n}$ we have $\|
A \|_{HS} = \sqrt{\sum_{i,j=1}^n |A_{i,j}|^2}$.

\begin{lemma}
Let $2 \leq k \leq n$ be integers, let $\alpha \geq 0$, and assume
that $0 < \lambda \leq 1$ is such that $k = n^{\lambda}$. Suppose
that $f: \RR^n \rightarrow [0, \infty)$ is an isotropic,
log-concave function and define $g = f * \gamma_{n}[n^{-\alpha
\lambda} ]$.
 Let $E_0 \subseteq \RR^n$ be a
$k$-dimensional subspace, and let $x_0 \in E_0$ be a point with
$|x_0| \leq 10 \sqrt{k}$. Then, for any   $U_1, U_2 \in SO(n)$,
$$ \left| M_{g, E_0, x_0}(U_1) - M_{g, E_0, x_0}(U_2)  \right| \leq
C n^{\lambda (2 \alpha + 2)}  \cdot d(U_1, U_2),
$$ where $C > 0$ is a universal constant.
\label{lemma_lipshitz}
\end{lemma}

\emph{Proof:} We abbreviate $M(U) = M_{g, E_0, x_0}(U)$ for $U \in
SO(n)$. We need to show that $M$ is $C n^{\lambda (2 \alpha +
2)}$-Lipshitz on $SO(n)$. By rotational-symmetry (see
(\ref{eq_642})), it is enough to show that
\begin{equation}
 |\nabla M(Id)| \leq C n^{\lambda (2 \alpha + 2)},
 \label{eq_549}
\end{equation}
where $\nabla M(Id)$ is the Riemannian gradient of $M: SO(n)
\rightarrow \RR$ at $Id$, and $|\nabla M(Id)|$ is its  length.
 Fix an orthonormal basis $e_1,...,e_k \in E_0$. For $v =
(v_1,...,v_k) \in (\RR^n)^k$, let $A_v: E_0 \rightarrow \RR^n$
stand for the unique linear operator with $A_v(e_i) = v_i$ for
$i=1,...,k$. Define,
$$ G(v_1,...,v_k) = \log \int_{A_vx_0 + (A_vE_0)^{\perp}} g(x) dx $$
which is well-defined because $x_0 \in E_0$.  Note that by
(\ref{eq_1141}), for any $U \in SO(n)$,
$$ M(U) = \log \int_{U x_0 + (U E_0)^{\perp}} g(x) dx
= G(U e_1,...,U e_k).
$$
Furthermore, for any $U \in SO(n)$,
$$ \sum_{i=1}^k |U e_i - e_i|^2  \leq \| U - Id \|_{HS}^2
\leq d(U, Id)^2.
$$
Therefore, to prove (\ref{eq_549}), it is sufficient to
demonstrate that
\begin{equation}
 |\nabla G(e_1,...,e_k)| \leq C n^{\lambda (2 \alpha + 2)}
 = C k^{2 \alpha + 2},
\label{improve1}
\end{equation}
where $\nabla G$ is  the usual gradient of the function $G:
(\RR^n)^k \rightarrow \RR$ in the Euclidean space $(\RR^n)^k$.
 For $i=1,...,k$ and $v \in \RR^n$ we set
$F_i(v) = G(e_1,...,e_{i-1}, v, e_{i+1}, ..., e_k)$. Then
$$
 |\nabla G(e_1,..,e_k)|^2 = \sum_{i=1}^k |\nabla F_i(e_i)|^2
 $$
(note that $G$ is a smooth function in a neighborhood of
$(e_1,...,e_k) \in (\RR^n)^k$, since $g$ is $C^{\infty}$-smooth
and decays exponentially fast at infinity). Therefore, it is
sufficient to prove that for any $i=1,..,k$,
\begin{equation} |\nabla F_i(e_i)| \leq C k^{2 \alpha + 3/2}.
\label{improve2}
 \end{equation}
By symmetry, it is enough to focus on the case $i = 1$. Thus, we
denote $F(v) = F_1(v) = G(v, e_2,...,e_k)$. Then $F: \RR^n
\rightarrow \RR$ is a smooth function in a neighborhood of $e_1$,
and our goal is to prove that
$$ |\nabla F(e_1)| \leq C k^{2 \alpha + 3/2}. $$
Equivalently, fix an arbitrary $v \in \RR^n$ with $|v| = 1$. To
prove the lemma, it suffices to show that
\begin{equation}
\left. \frac{d}{dt} F(e_1 + t v) \right|_{t =0} = \langle \nabla
F(e_1), v \rangle \leq C k^{2 \alpha + 3/2}. \label{eq_632}
\end{equation}
We thus focus on proving (\ref{eq_632}). Denote $E = sp (E_0 \cup
\{ v \})$ where $sp$ denotes linear span, and let $\bar{g} =
\pi_E(g)$. For $t \in \RR$ and $x \in E_0$ set $A_t(x) = x + t
\langle x, e_1 \rangle v$. Then for any $|t| \leq 1/2$ we have
 $A_t(E_0) \subseteq E$ and
$$ F(e_1 + t v) = \log \pi_{A_t(E_0)}(g)  (A_t x_0) = \log \pi_{A_t(E_0)}
\left( \bar{g} \right) (A_t x_0). $$ We have thus reduced our
$n$-dimensional problem to a $(k+1)$-dimensional problem; the
function $\bar{f} := \pi_E(f)$ is isotropic and log-concave on $E$
(by Pr\'ekopa-Leindler), and
\begin{equation}
 \bar{g} = \pi_E(g) = \pi_E \left( f * \gamma_n[n^{-\alpha
\lambda}] \right) = \bar{f} * \gamma_{E}
[k^{-\alpha}].
\label{eq_1138}
\end{equation}
We divide the proof of (\ref{eq_632}) into two cases. Suppose
first that $v \in E_0$. In this case, $E = E_0 = A_t(E_0)$ for all
$|t| < 1/2$. Therefore, for $|t| < 1/2$,
$$ F(e_1 + t v) = \log \bar{g}(A_t x_0) = \log \bar{g}(x_0 + t \langle x_0, e_1 \rangle
v).
$$
The subspace $E$ is $k$-dimensional. We may apply Lemma
\ref{lem_349} for $k, \bar{f}, \bar{g}, x_0$ and $v = k^{-\alpha}$
because of (\ref{eq_1138}). By the conclusion of that lemma,
$$ \left. \frac{d}{dt} F(e_1 + t v) \right|_{t =0} = \left \langle
\nabla \log \bar{g}(x_0), v \right \rangle \cdot \langle x_0, e_1
\rangle \leq C k^{ \alpha + 1/2} |x_0| \leq C^{\prime} k^{\alpha +
1},
$$
since $|v| = 1$ and $|x_0| \leq 10 \sqrt{k}$. Thus (\ref{eq_632})
is proven for the case where $v \in E_0$.

\medskip From this point on and until the end of the proof,
we suppose that $v \not \in E_0$ and our aim is to prove
(\ref{eq_632}).
 Then $\dim(E)
= k+1$, and for $|t| < 1/2$,
$$
 F(e_1 + t v) = \log \pi_{A_t(E_0)}
(\bar{g}) (A_t x_0) = \log \int_{-\infty}^{\infty} \bar{g}
\left(A_t x_0 + r \theta_t \right) dr,
$$
where $\theta_t$ is a unit vector that is orthogonal to the
hyperplane $A_t(E_0)$ in $E$. There exist two such unit vectors,
and we may select any of them. For concreteness, we choose
$$ \theta_t = \frac{ t |Proj_{E_0^{\perp}} v|^2 e_1 - (1 + t
\langle v, e_1 \rangle) Proj_{E_0^{\perp}} v}{ | \, t
|Proj_{E_0^{\perp}} v|^2 e_1 - (1 + t \langle v, e_1 \rangle)
Proj_{E_0^{\perp}} v \, | }.
$$
An elementary geometric argument (or, alternately, a tedious
computation) shows that $\left. d \theta_t / dt \right|_{t=0} $ is
a vector whose length is at most one, since $|v| \leq 1$.

\medskip For $s, t \in \RR$ with $|s|, |t| \leq 1/2$, we define
$$  \bar{F}(s, t) = \log \int_{-\infty}^{\infty} \bar{g} \left(A_t x_0 + r \theta_s \right) dr, $$
which is a smooth function in $s$ and $t$. Then $\bar{F}(t,t) =
F(e_1 + t v)$. Recall that $\bar{g}$ is the convolution of an
isotropic, log-concave function with $\gamma_E[k^{-\alpha}]$,
according to (\ref{eq_1138}). Lemma \ref{rotate} implies that
$$  \left. \frac{\partial \bar{F}(s, t)}{\partial s} \right|_{s=t=0}
\leq C (k+1)^{3/2}  k^{2 \alpha}  \, \cdot \, \left| \, \left.
\partial \theta_s / \partial s \right|_{s=0} \, \right| \leq
C^{\prime} k^{2 \alpha + 3/2}.$$ Next, we estimate $\partial
\bar{F} /
\partial t$. Note that for any $x \in \RR^n$,
$$ \log \int_{-\infty}^{\infty} \bar{g}(x + r \theta_0) dr = \log
\int_{-\infty}^{\infty} \bar{g}(Proj_{E_0}(x) + r \theta_0) dr =
\log \pi_{E_0}(\bar{g})(Proj_{E_0}(x)). $$ Therefore, $\bar{F}(0,
t) = \log \pi_{E_0}(\bar{g}) (x_0 + t \langle x_0, e_1 \rangle
Proj_{E_0} (v))$ for any $|t| < 1/2$. The function
$\pi_{E_0}(\bar{g})$ is the convolution of an isotropic,
log-concave function with $\gamma_{E_0}[k^{-\alpha}]$. We appeal
to Lemma \ref{lem_349} and conclude that
$$  \left. \frac{\partial \bar{F}(s,t)}{\partial t} \right|_{t=s=0} \leq C k^{\alpha + 1/2}
 \cdot |\langle x_0, e_1 \rangle| \cdot
|Proj_{E_0} (v)| \leq C^{\prime} k^{\alpha + 1},
$$
since $|v| \leq 1$ and $|x_0| \leq 10 \sqrt{k}$. We have thus
shown that
$$ \left.  \frac{d F(e_1 + t v)}{dt} \right|_{t=0} =
\left. \frac{\partial \bar{F}(s,t)}{\partial s} \right|_{t=s=0} +
\left. \frac{\partial \bar{F}(s,t)}{\partial t} \right|_{t=s=0}
\leq C k^{2 \alpha + 3/2}. $$ This completes the proof of
(\ref{eq_632}) in the case where $v \not \in E_0$. Thus
(\ref{eq_632}) holds in all cases. The lemma is proven. \hfill
$\square$

\medskip The group $SO(n)$ admits a unique
Haar probability measure $\mu_n$, which is invariant under both
left and right translations. Our next proposition is a
concentration of measure inequality on the orthogonal group from
Gromov and Milman \cite{GM}, see also \cite[Section 6 and Appendix
V]{MS}. This measure-concentration inequality is deduced in
\cite{GM} from a very general isoperimetric inequality due to
Gromov, which requires only lower bounds on the Ricci curvature of
the manifold in question. In the specific case of the orthogonal
group, a more elementary proof of Proposition \ref{prop_651} may
be obtained by using two-point symmetrization (see
\cite{benyamini}).

\begin{proposition} Let $n \geq 1, \eps > 0, L > 0$ and let
$f: SO(n) \rightarrow \RR$ be such that
$$ f(U) - f(V) \leq L d(U, V) $$ for all $U, V \in SO(n)$. Denote
$M = \int_{SO(n)} f(U) d\mu_n(U)$. Then,
$$ \mu_n \left \{ U \in SO(n) ; \left| f(U) - M \right| \geq \eps \right \}
\leq C \exp \left(- c n \eps^2 / L^2 \right), $$ where $C,c > 0$
are universal constants. \label{prop_651}
\end{proposition}

Milman's principle from \cite{mil_dvo} states, very roughly, that
Lipshitz functions on certain high-dimensional mathematical
structures are approximately constant when restricted to typical
sub-structures. Behind this principle there usually stands a
concentration of measure inequality, such as Proposition
\ref{prop_651}. Our next lemma is yet another manifestation of
Milman's principle, whose proof is rather similar to the original
argument in \cite{mil_dvo}.

\medskip Recall that $G_{n,k}$ stands for the grassmannian of all
$k$-dimensional subspaces in $\RR^n$. There is a unique rotationally-invariant
probability measure on $G_{n,k}$, denoted by $\sigma_{n,k}$.
Whenever we say that $E$ is a random $k$-dimensional subspace
in $\RR^n$, and whenever we say that $U$ is a random rotation
in $SO(n)$, we relate to the probability measures $\sigma_{n,k}$
and $\mu_n$, respectively.
For a subspace $E \in G_{n,k}$, let $S(E) = \{ x \in E ;|x| = 1\}$
stand for the unit sphere in $E$, and let $\sigma_E$ denote the
unique rotationally-invariant probability measure on $S(E)$.

\begin{lemma}
Let $2 \leq k \leq n$ be integers, let $0 \leq \alpha \leq 10^{5},
-10 \leq \eta \leq 10^5, 0 < u < 1$, and assume that
\begin{equation}
0 < \lambda \leq \min \left \{  1 / (4 \alpha + 2 \eta + 5.01)  ,
u / (4 \alpha + 2 \eta + 4) \right \} \label{lambda}
\end{equation}
is such that $k = n^{\lambda}$. Suppose that $f: \RR^n \rightarrow
[0, \infty)$ is an isotropic, log-concave function and define $g =
f * \gamma_{n}[n^{-\lambda \alpha} ]$.

\smallskip Let $E \in G_{n,k}$ be a random subspace. Then, with
probability greater than $1 - C e^{-c n^{1 - u}}$,
$$ \left| \log \pi_E(g)(x_1) - \log \pi_E(g)(x_2) \right| \leq \frac{1}{k^{\eta}} $$
for all $x_1, x_2 \in E$ with $|x_1| = |x_2| \leq 10 \sqrt{k}$.
Here, $c, C
> 0$ are universal constants. \label{lem_628}
\end{lemma}

\emph{Proof:} Fix a $k$-dimensional subspace $E_0 \subseteq
\RR^n$. By rotational-invariance, for any $r > 0$ and a unit
vector $v \in E_0$,
\begin{eqnarray}
\label{eq_451} \lefteqn{ \int_{SO(n)} M_{g, E_0, r v}(U) d
\mu_n(U) = \int_{SO(n)} \log \pi_{U(E_0)}(g)  (r U v) d \mu_n(U) }
\\ & = & \int_{G_{n,k}} \int_{S(E)} \log \pi_E(g)(r \theta) d
\sigma_E(\theta) d \sigma_{n,k}(E). \phantom{aaaaaaaaaaaaaaaaa}
\nonumber
\end{eqnarray}
Consequently, for $r > 0$ we may define
\begin{equation}
 M(r) = \int_{SO(n)} M_{g, E_0, r v}(U) d \mu_n(U),
 \label{eq_204}
\end{equation}
  where
$v \in E_0$ is an arbitrary unit vector, and the definition does
not depend on the choice of the unit vector $v \in E_0$.

\smallskip According to Lemma \ref{lemma_lipshitz}, for any $x \in
E_0$ with $|x| \leq 10 \sqrt{k}$, the function
$$ U \mapsto M_{g, E_0, x}(U) $$
is $C k^{2  \alpha + 2}$-Lipshitz on $SO(n)$. Therefore, by
Proposition \ref{prop_651} and by (\ref{eq_204}), for any $x \in
E_0$ with $|x| \leq 10 \sqrt{k}$ and for any $\eps > 0$,
\begin{equation}
 \mu_n \left \{ U \in SO(n) ; \left| M_{g, E_0, x}(U) - M(|x|) \right| > \eps
\right \} \leq C \exp \left(-c \frac{\eps^2 n}{k^{4 \alpha + 4}}
\right). \label{eq_206}
\end{equation}
Let $\eps, \delta
> 0$ be some small numbers to be specified later. Let $\cN  \subset
10 \sqrt{k} D^n
 \cap E_0$ be an $\eps$-net for $10 \sqrt{k} D^n \cap E_0$ of at most $(C \sqrt{k} / \eps)^k$
elements
 (see, e.g. \cite[Lemma 4.10]{pisier}). That
is, for any $x \in 10 \sqrt{k} D^n \cap E_0$ there exists $y \in
\cN$ with $|x - y| \leq \eps$.
 Suppose $U \in SO(n)$
is such that
\begin{equation} \left| \, M_{g, E_0, x}(U) - M (|x|) \, \right| \leq \delta \ \
\ \text{for all} \ x \in \cN. \label{eq_332}
\end{equation}
Denote $E = U(E_0)$, and let $\cN^{\prime} = \{ U x ; x \in \cN
\}$. Then (\ref{eq_332}) and the definition (\ref{eq_1141}) imply
that
\begin{equation}
 \left| \, \log \pi_E(g)(x) - M(|x|) \, \right| \leq \delta
\ \ \ \text{for all} \ x \in \cN^{\prime}. \label{eq_345}
\end{equation}
 The function
$\pi_E(f)$ is isotropic and log-concave, and $\pi_E(g) = \pi_E(f)
* \gamma_E[k^{-\alpha}]$. We may thus apply Lemma \ref{lem_349}
and conclude that the function $x \mapsto \log \pi_E(f)(x)$ is $C
k^{\alpha + 1/2}$-Lipshitz on $10 \sqrt{k} D^n \cap E$. Therefore,
from (\ref{eq_451}) and (\ref{eq_204}) we deduce that $r \mapsto
M(r)$ is $C k^{\alpha + 1/2}$-Lipshitz on the interval $(0, 10
\sqrt{k})$. Since $\cN^{\prime}$ is an $\eps$-net for $10 \sqrt{k}
D^n \cap E$, we infer from (\ref{eq_345}) that
$$
 \left| \, \log \pi_E(g)(x) - M(|x|) \, \right| \leq \delta + 2 \eps
 C k^{\alpha + 1/2} \ \ \ \text{for all} \ x \in 10 \sqrt{k} D^n \cap E. $$
To summarize, if $U \in SO(n) $ is such that (\ref{eq_332}) holds,
then for all $x, x^{\prime} \in E = U(E_0)$,
\begin{equation}
\left| \, \log \pi_E(g)(x) - \log \pi_E(g)(x^{\prime}) \, \right|
\leq 2\delta + 4 \eps
 C k^{ \alpha + 1/2} \ \ \ \text{when} \ |x| = |x^{\prime}| \leq 10
 \sqrt{k}.
 \label{eq_346}
\end{equation}
Recall the estimate (\ref{eq_206}). Recall that $\cN$ is a subset
of $10 \sqrt{k} D^n \cap E_0$ of cardinality at most $(C^{\prime}
\sqrt{k} / \eps)^k$. Therefore, the probability of a random
rotation $U \in SO(n)$ to satisfy (\ref{eq_332}) is greater than
$$
 1 - (C^{\prime} \sqrt{k} / \eps)^k  \cdot C \exp \left(-c^{\prime}
\delta^2 k^{-4 \alpha - 4} n \right). $$ Set $\delta = k^{-\eta} /
10 $ and $\eps = \delta k^{-\alpha-1/2} / C$ where $C$ is the
constant from (\ref{eq_346}). Since $$ k \leq \min \left \{ n^{1 /
(4 \alpha + 2 \eta + 5.01)}, n^{u / (4 \alpha + 2 \eta + 4 )}
\right \}
$$ by (\ref{lambda}), then $\delta^2 k^{-4 \alpha - 4} n > c n^{1 - u}$ and
also
$$ 1 - (C^{\prime} \sqrt{k} / \eps)^k  \cdot \exp
\left(-c^{\prime} \delta^2 k^{-4 \alpha - 4} n \right) \geq 1 -
\bar{C} e^{-\bar{c} n^{1 - u}}. $$ We conclude that if $U$ is a
random rotation in $SO(n)$, then (\ref{eq_332}) holds  with
probability greater than $1 - \bar{C} \exp(-\bar{c} n^{1 - u})$.
Whenever $U \in SO(n)$ satisfies (\ref{eq_332}), the subspace $E =
U(E_0)$ necessarily satisfies (\ref{eq_346}). Hence, with
probability greater than $1 - \bar{C} \exp(-\bar{c} n^{1 - u})$ of
selecting $U \in SO(n)$, for any $x, x^{\prime} \in U(E_0)$,
$$ \left| \, \log \pi_{U(E_0)}(g)(x) - \log \pi_{U(E_0)}(g)(x^{\prime}) \, \right|
\leq \frac{1}{k^{\eta}}
 \ \ \ \text{when} \ |x| = |x^{\prime}| \leq 10 \sqrt{k}. $$
Note that the subspace $U(E_0)$ is distributed uniformly on
$G_{n,k}$. The proof is complete.\hfill $\square$

\section{Almost-radial marginals}

Suppose that $f: \RR^n \rightarrow [0, \infty)$ is an isotropic,
log-concave function. Lemma \ref{lem_628} states that a typical
$n^c$-dimensional marginal of $f$
 is approximately radial, after
convolving with a gaussian. In this section we will show -- mostly
by referring to \cite{CLT} -- that a large portion of the mass of
this typical marginal is located in a very thin spherical shell.

\smallskip Let $f: [0, \infty) \rightarrow [0, \infty)$ be a log-concave
function with $0 < \int_0^{\infty} f < \infty$ that is continuous
on $[0, \infty)$ and $C^2$-smooth on $(0, \infty)$. As in
\cite{CLT}, for $p
> 1$, we denote by $t_p(f)$ the unique $t > 0$ for which $f(t) > 0$
and also $$ (\log  f)^{\prime}(t) = \frac{f^{\prime}(t)}{f(t)} =
-\frac{p-1}{t}. $$ The quantity  $t_p(f)$ is well-defined
according to \cite[Lemma 4.3]{CLT}.  The following lemma asserts
that most of the mass of $t \mapsto t^{p-1} f(t)$ is located in a
small neighborhood of $t_p(f)$. We refer the reader to \cite[Lemma
4.5]{CLT} for the proof.

\begin{lemma} Let $p \geq 2$, and let $f: [0, \infty) \rightarrow [0, \infty)$
be a continuous, log-concave function, $C^2$-smooth on $(0,
\infty)$, with $0 < \int_0^{\infty} f < \infty$. Then for all $0 \leq
\eps \leq 1$,
$$
\int_{t_p(f) (1 - \eps)}^{t_p(f) (1 + \eps)} t^{p-1} f(t) dt \geq
\left(1 - C e^{-c \eps^2 p} \right) \int_{0}^{\infty} t^{p-1} f(t)
dt, $$
 where $C, c > 0$ are universal constants. \label{lem_1153}
\end{lemma}

Next, we analyze log-concave densities that are almost-radial in
the sense of Lemma \ref{lem_628}.

\begin{lemma} Let $n \geq C$
and let $f: \RR^n \rightarrow [0, \infty)$ be a $C^{2}$-smooth,
log-concave probability density. Let $X$ be a random vector in
$\RR^n$ whose density is $f$, and assume that $\EE X = 0$ and $n
\leq \EE |X|^2 \leq 2n$. Suppose that $\delta > 0$ is such that
\begin{equation}
 \left| \log f(x_1) - \log f(x_2) \right| \leq
\delta n \label{eq_838}
\end{equation}
 for all $x_1, x_2 \in \RR^n$ with $|x_1| = |x_2| \leq
10 \sqrt{n}$. Denote $r = \sqrt{\EE |X|^2}$. Then for all $\eps
> 0$ with $C \sqrt{\delta} \leq \eps \leq 1$,
$$ Prob \left \{ \left|  \frac{|X|}{r} - 1   \right| > \eps \right
\} \leq C e^{-c \eps^2 n }. $$
 Here, $C, c > 0$ are universal constants. \label{lem_853}
\end{lemma}

\emph{Proof:} We may assume that $\delta < 10^{-3}$; otherwise,
there is no $\eps > 0$ with $C \sqrt{\delta} \leq \eps \leq 1$ for
a sufficiently large universal constant $C$. For $\theta \in
S^{n-1}$ and $r \geq 0$ we denote $f_{\theta}(r) = f(r \theta)$.
Since $\int f =1$ then $f$ decays exponentially fast at infinity
(e.g., \cite[Lemma 2.1]{psitwo}). Consequently,
$t_{n}(f_{\theta})$ is well-defined for all $\theta \in S^{n-1}$.
Let $M > 0$ be such that
\begin{equation}
 \int_{ \{ x \in \RR^n ; |x| \leq M \} } f(x) dx =
\frac{2}{3} \int_{\RR^n} f(x) dx = \frac{2}{3}. \label{eq_1211}
\end{equation}
Since $\int_{\RR^n} |x|^2 f(x) dx \geq n$, then Borell's lemma
(e.g. \cite[Appendix III.3]{MS}) implies that $M \geq \sqrt{n} /
10$. Additionally, since $\int_{\RR^n} |x|^2 f(x) dx \leq 2 n$,
then necessarily $M \leq 3 \sqrt{n}$, by Markov's inequality. We
integrate (\ref{eq_1211}) in polar coordinates and obtain
\begin{equation} \int_{S^{n-1}} \int_0^{M} f_{\theta}(t)
t^{n-1} dt d \theta = \frac{2}{3} \int_{S^{n-1}} \int_0^{\infty}
f_{\theta}(t) t^{n-1} dt d \theta. \label{eq_1212}
\end{equation} We claim that there exists $\theta_0 \in S^{n-1}$
with
\begin{equation}
\frac{1}{20} \leq  \frac{t_n(f_{\theta_0})}{\sqrt{n}} \leq 6.
\label{eq_204_}
\end{equation}
Otherwise, by continuity, either $t_n(\theta) > 6\sqrt{n}$ for all
$\theta \in S^{n-1}$ or else $t_n(\theta) < \sqrt{n} / 20$ for all
$\theta \in S^{n-1}$. In the first case, for all $\theta \in
S^{n-1}$ we have $t_n(\theta) > 6 \sqrt{n} \geq 2 M$, and by Lemma
\ref{lem_1153},
$$ \forall \theta \in S^{n-1}, \ \ \ \ \ \int_0^{M} f_{\theta}(t) t^{n-1} dt \leq
\int_0^{t_n(\theta) / 2} f_{\theta}(t) t^{n-1} dt < \frac{2}{3}
\int_0^{\infty} f_{\theta}(t) t^{n-1} dt, $$ provided that $n >
C$, in contradiction to (\ref{eq_1212}). Similarly, in the second
case,
 for all $\theta \in
S^{n-1}$, we have $2 t_n(\theta) < \sqrt{n} / 10 \leq M$ and by Lemma \ref{lem_1153},
$$ \forall \theta \in S^{n-1}, \ \ \ \ \ \int_0^{M} f_{\theta}(t) t^{n-1} dt \geq
\int_0^{2 t_n(\theta)} f_{\theta}(t) t^{n-1} dt > \frac{2}{3}
\int_0^{\infty} f_{\theta}(t) t^{n-1} dt, $$ in contradiction to
(\ref{eq_1212}). We have thus proven that there exists $\theta_0
\in S^{n-1}$ such that (\ref{eq_204_}) holds. Fix such $\theta_0
\in S^{n-1}$. Denote $\vphi_0(t) = \log f_{\theta_0}(t)$ for $t
\geq 0$ (where $\log 0 = -\infty$) and $r_0 = t_n(f_{\theta_0})$.
Then,
\begin{equation}
\sqrt{n} / 20 \leq r_0 \leq 6 \sqrt{n}. \label{eq_702}
\end{equation}
Fix $\theta \in S^{n-1}$ and denote
$r = t_n(f_{\theta})$
and  $\vphi(t) = \log
f_{\theta}(t)$ for $t \geq 0$.
Then
$\vphi^{\prime}(r) = -(n-1) / r$. We will prove that
\begin{equation}
r_0 / r \leq 1+ 60 \sqrt{\delta}. \label{eq_818}
\end{equation}
Indeed, assume the contrary. Then $r_0 > r (1+ 60 \sqrt{\delta})$.
Since $r_0 \geq \sqrt{n} / 20$ and $\sqrt{\delta} < 1 / 30$, then
necessarily \begin{equation} r_0 - r
> \sqrt{\delta n}. \label{eq_936} \end{equation} Recall that $\vphi$ and $\vphi_0$ are concave functions,
hence their derivatives are non-increasing. Therefore for all $t
\in [r, r_0]$,
\begin{eqnarray}
\label{eq_1247} \lefteqn{ \vphi^{\prime}(t) \leq \vphi^{\prime}(r)
=
-\frac{n-1}{r} < -\frac{n-1}{r_0} \left (1 + 60 \sqrt{\delta} \right) }\\
& \leq & -\frac{n-1}{r_0} - 6 \sqrt{\delta n} =
\vphi_0^{\prime}(r_0) - 6 \sqrt{\delta n} \leq \vphi_0^{\prime}(t)
- 6 \sqrt{\delta n}, \nonumber
\end{eqnarray}
where we used the fact that $(n-1) / r_0 \geq (n-1) / (6 \sqrt{n})
\geq \sqrt{n} / 10$ by (\ref{eq_702}). Note that $r < r_0 \leq 6
\sqrt{n}$, and hence (\ref{eq_838}) implies that $|\vphi_0(t) -
\vphi(t)| \leq \delta n$ for all $t \in [r, r_0]$. However, by
(\ref{eq_1247}) and (\ref{eq_936}),
$$ [\vphi_0(r_0) - \vphi(r_0)] - [\vphi_0(r) - \vphi(r)] = \int_{r}^{r_0} \vphi_0^{\prime}(t) - \vphi^{\prime}(t) dt
> 6 \sqrt{\delta n} (r_0 - r) >  6 \delta n, $$ in contradiction to (\ref{eq_838}).
 Thus (\ref{eq_818}) is proven. Next we will demonstrate
that
\begin{equation}
r / r_0 \leq 1+ 200 \sqrt{\delta}. \label{eq_208}
\end{equation}
The proof of (\ref{eq_208}) is very similar to the proof  of
(\ref{eq_818}). Assume on the contrary that (\ref{eq_208}) does
not hold. Then $r > r_0 + 10 \sqrt{\delta n}$. Denote $\bar{r} =
r_0 + 10 \sqrt{\delta n}$. Since $\vphi^{\prime}$ is
non-increasing, then $\vphi^{\prime}(\bar{r}) \geq
\vphi^{\prime}(r) = -(n-1) / r > -(n-1) / \bar{r}$. Hence, for $t
\in [r_0, \bar{r}]$,
$$
\vphi_0^{\prime}(t) \leq \vphi_0^{\prime}(r_0) = -\frac{n-1}{r_0}
< -\frac{n-1}{\bar{r}} - \frac{\sqrt{\delta n}}{5} <
\vphi^{\prime}(\bar{r}) - \frac{\sqrt{\delta n}}{5} \leq
\vphi^{\prime}(t) - \frac{\sqrt{\delta n}}{5}, $$
 according to
(\ref{eq_702}). As before, this entails
$$ [\vphi(\bar{r}) - \vphi_0(\bar{r})] - [\vphi(r_0) -
\vphi_0(r_0)] > (\sqrt{\delta n} / 5) \cdot (\bar{r} - r_0) = 2
\delta n,
$$ in contradiction to (\ref{eq_838}).
 To summarize, in (\ref{eq_818}) and (\ref{eq_208})
we proved that for all $\theta \in S^{n-1}$,
\begin{equation}
\left( 1 -  C \sqrt{\delta} \right) r_0 \leq  t_n(f_{\theta}) \leq
\left( 1+  C \sqrt{\delta} \right) r_0. \label{eq_1048_}
\end{equation} We may assume that $10 C \sqrt{\delta} < 1$, where $C$
is the constant from (\ref{eq_1048_}). Let $\eps
> 0$ satisfy $10 C \sqrt{\delta} < \eps < 1$. According to (\ref{eq_1048_}), for any $\theta \in S^{n-1}$,
$$ r_0(1-\eps) \leq t_n(f_{\theta})(1-\eps/2)
\ \ \ \ \ \text{and} \ \ \ \ \ r_0(1+\eps) \geq
t_n(f_{\theta})(1+\eps/2).
$$
 Integration in polar coordinates yields
\begin{eqnarray*}
\lefteqn{ Prob \left \{ \left| \frac{|X|}{r_0} - 1  \right| \leq
\eps \right \} = \int_{S^{n-1}} \int_{r_0(1-\eps)}^{r_0(1+\eps)}
f(t \theta) t^{n-1} dt d\theta } \\ & \geq & \int_{S^{n-1}}
\int_{t_n(f_{\theta})(1-\eps/2)}^{t_n(f_{\theta})(1+\eps/2)}
f_{\theta}(t) t^{n-1} dt d\theta \\ & \geq & \left(1 - C^{\prime}
e^{-c^{\prime} \eps^2 n} \right) \int_{S^{n-1}} \int_0^{\infty}
f_{\theta}(t) t^{n-1} dt d\theta
 =  1 - C^{\prime} e^{-c^{\prime} \eps^2 n},
\end{eqnarray*}
where we used  Lemma \ref{lem_1153}. Therefore, when $\tilde{C}
\sqrt{\delta} < \eps < 1$, \begin{equation}  Prob \left \{ \left|
\frac{|X|^2}{r_0^2} - 1  \right| \geq \eps \right \} \leq C e^{-c
\eps^2 n}. \label{eq_933} \end{equation} As in the proof of Lemma
4.6 in \cite{CLT}, we use (\ref{eq_933}) and Borell's lemma to
obtain
\begin{eqnarray} \label{eq_946} \lefteqn{\left|  \frac{\EE |X|^2}{r_0^2} - 1  \right|
\leq \EE \left|  \frac{|X|^2}{r_0^2} - 1  \right| =
\int_0^{\infty} Prob \left \{ \left| \frac{|X|^2}{r_0^2} - 1
\right| \geq t \right \} dt } \\ & \leq & \tilde{C} \sqrt{\delta}
+ \int_{\tilde{C} \sqrt{\delta}}^{1} C e^{-c t^2 n} dt +
C^{\prime} \int_1^{\infty} \min \{ e^{-c n}, e^{-c^{\prime}
\sqrt{t}} \} dt \leq \bar{C} \sqrt{\delta + \frac{1}{n}}.
\nonumber
\end{eqnarray}
Recall that we denote $r = \sqrt{\EE |X|^2}$. From (\ref{eq_933})
and (\ref{eq_946}) we conclude that
$$ Prob \left \{ \left|
\frac{|X|^2}{r^2} - 1  \right| \geq \eps \right \} \leq C e^{-c
\eps^2 n} $$ for all $\eps > 0$ with $\tilde{C} \sqrt{\delta} <
\eps < 1$. This completes the proof. \hfill $\square$

\begin{lemma} Let $n \geq 2$ be an integer, let
\begin{equation}
 \max \left \{ \frac{1}{10}, \frac{1}{\log n} \right \} \leq
\lambda \leq \frac{1}{3.01} - \frac{C^{\prime}}{\log n},
\label{eq_1048}
\end{equation}
 and assume that
$X$ is a random vector in $\RR^n$ with an isotropic, log-concave
density. Let $Y$ be a standard gaussian random vector in $\RR^n$,
independent of $X$.
 Then
$$
Prob \left \{ \left| \frac{|X + Y|}{\sqrt{2n}} - 1  \right| > C
n^{(3.01 \lambda - 1) / 4 } \right \} \leq C \exp \left( -c
n^{(5.01 \lambda - 1) / 2} \right).
$$
Here, $c, C, C^{\prime} > 0$ are universal constants.
\label{lem_634}
\end{lemma}

\emph{Proof:} We may clearly assume that $n \geq \tilde{C}$.
Denote by $f$ the density of the random vector $X$. Then $f: \RR^n
\rightarrow [0, \infty)$ is isotropic and log-concave. Define $g =
f * \gamma_{n}[1]$. Let $k$ be the maximal integer such that $k
\leq n^{\lambda}$. Then $k \geq \max \{2, n^{\lambda} / 2 \}$
because of (\ref{eq_1048}). Define
 $$ \eta = \frac{1}{2 \lambda} - 2.505. $$
 Then $-1 \leq \eta \leq 10$.
We apply Lemma \ref{lem_628} for $\alpha = 0$, for  $\eta$ as was
just defined and for $u = 1 - 1.01 \lambda$. Note that
$$ \lambda \leq \min \{  1 / (2 \eta + 5.01), u / (2 \eta + 4) \}.
$$
Thus the appeal to Lemma \ref{lem_628} is legitimate. By the
conclusion of Lemma \ref{lem_628}, there exists $\cE \subseteq
G_{n,k}$ with
$$ \sigma_{n,k} (\cE) >  1 - C e^{-c n^{1-u}}  \geq 1 - C e^{-c n^{\lambda}} $$
and with the following property: For any $E \in \cE$ and $x_1, x_2
\in E$,
$$
 \left| \log \pi_E(g)(x_1) - \log \pi_E(g)(x_2) \right| \leq
k^{-\eta} \ \ \ \ \text{when} \ \  |x_1| = |x_2| \leq 10 \sqrt{k}.
$$
Equivalently, denote $\delta = k^{-\eta - 1}$. For an appropriate
choice of a large universal constant $C^{\prime}$ in
(\ref{eq_1048}), we have that $C \sqrt{\delta} \leq 1$ where $C$
is the constant from Lemma \ref{lem_853}. Then for any $E \in \cE$
and $x_1, x_2 \in E$,
\begin{equation}
 \left| \log \pi_E(g)(x_1) - \log \pi_E(g)(x_2) \right| \leq
\delta k \ \ \ \ \text{when} \ \  |x_1| = |x_2| \leq 10 \sqrt{k}.
 \label{eq_356}
\end{equation}
Fix a subspace $E \in \cE$. The function $\pi_E(g)$ is
$C^2$-smooth (it is a convolution with a gaussian), log-concave
(by Pr\'ekopa-Leindler), and it satisfies (\ref{eq_356}). The
random vector  $Proj_E(X + Y)$ is distributed according to the
density $\pi_E(g) = \pi_E(f) * \gamma_E[1]$ in the subspace $E$.
Furthermore,
$$ \EE Proj_E(X + Y) =
0 \ \ \ \ \ \ \text{and} \ \ \ \ \ \  \EE |Proj_E(X + Y)|^2 =  2
k.
$$
We have thus verified the assumptions of Lemma \ref{lem_853} for
the function $\pi_E(g)$, the random vector $Proj_E(X + Y)$ and the
number $\delta$. By the conclusion of Lemma \ref{lem_853} (the
case $\eps = C \sqrt{\delta} \leq 1$),
\begin{equation} Prob \left \{ \left| \frac{|Proj_E(X  + Y)|}
{\sqrt{2k}} - 1  \right| > C \sqrt{\delta} \right \} \leq
C^{\prime} e^{-c^{\prime} \delta k }. \label{eq_358}
\end{equation}
The rest of the argument is similar to the proof of Theorem 1.4 in
\cite{CLT}; see the derivation involving formulas (49), (50) in
\cite[Section 4]{CLT}. We have proven that (\ref{eq_358}) holds
for all $E \in \cE$. Recall that $\sigma_{n,k}(\cE) > 1 - C
\exp(-c n^{\lambda} )$. Let $E \in G_{n,k}$ be a random subspace,
independent of $X$ and $Y$. Then
\begin{equation} Prob \left \{ \left| \frac{|Proj_E(X  + Y)|}
{\sqrt{2k}} - 1  \right| > C \sqrt{\delta} \right \} \leq
C^{\prime} e^{-c^{\prime} \delta k } + C e^{-c n^{\lambda}}.
\label{eq_1103_}
\end{equation}
However, according to the Johnson-Lindenstrauss dimension
reduction lemma (see, e.g. \cite[Lemma 4.8]{CLT}),
\begin{equation}
Prob \left \{   \left| \frac{ |Proj_E(X + Y)|}{|X + Y|} -
\sqrt{\frac{k}{n}} \right| \geq \sqrt{\delta} \cdot
\sqrt{\frac{k}{n}} \right \} \geq 1 - C e^{-c \delta k}.
\label{eq_1100}
\end{equation}
From (\ref{eq_1103_}) and (\ref{eq_1100}) we obtain
\begin{equation}
 Prob \left \{ \left| \frac{|X  + Y|} {\sqrt{2n}} - 1  \right|
> C \sqrt{\delta} \right \} < C e^{-c \delta k} + C e^{-c
n^{\lambda}} < C^{\prime} e^{-c \delta k}, \label{eq_1209}
\end{equation} since $\delta < 1$ and $k \leq n^{\lambda}$.
Recall that $ n^{\lambda} / 2 \leq k \leq n^{\lambda}$ and that
$\delta = k^{(3.01 - 1 / \lambda) / 2}$. The lemma  follows from
(\ref{eq_1209}). \hfill $\square$

\begin{theorem} Let $n \geq 1$ be an integer and let $X$ be
a random vector with an isotropic, log-concave density in $\RR^n$.
Then for all $0 \leq \eps \leq 1$,
\begin{equation}
 Prob \left \{ \left| \, \frac{|X|}{\sqrt{n}} -1 \, \right| \geq
\eps \right \} \leq C \exp \left( -c \eps^{3.33} n^{0.33} \right),
\label{eq_607}
\end{equation} where $c, C > 0$ are universal constants.
  \label{prop_408}
\end{theorem}

\emph{Proof:} We may assume that $n \geq C$. Let $Y$ be a standard
gaussian random vector in $\RR^n$, independent of $X$. Let
$n^{-0.1} \leq \eps \leq \bar{c}$, for a sufficiently small
universal constant $\bar{c} > 0$. Then the real number $\lambda$
defined by the equation $\eps = n^{(3.01 \lambda - 1) / 4}$
satisfies (\ref{eq_1048}). Consequently, by Lemma \ref{lem_634},
 $$ Prob \left \{
\left| \frac{|X + Y|}{\sqrt{2n}} - 1 \right|
> C \eps  \right \} \leq  C \exp \left( -c n^{(5.01 \lambda -
1) / 2} \right)  \leq  C^{\prime} \exp \left( -c^{\prime}
\eps^{3.33} n^{0.33} \right). $$ By adjusting the constants, we
conclude that
 for all $0 < \eps < 1$,
\begin{equation}
 Prob \left \{ \left| \frac{|X + Y|}{\sqrt{2n}} - 1 \right| > \eps
\right \} \leq \tilde{C} \exp \left( -\tilde{c} \eps^{3.33}
n^{0.33} \right). \label{eq_119} \end{equation} The random vector
$X$ has an isotropic, log-concave density. The standard gaussian
random vector $Y$ is independent of $X$. The simple argument that
leads from (\ref{eq_119}) to (\ref{eq_607}) was described in great
detail in the proof of Proposition 4.1 in \cite{CLT}; see the
derivation involving formulas (40),...,(45) in \cite[Section
4]{CLT}. We will not repeat that argument here. This completes the
proof. \hfill $\square$

\medskip
\emph{ Proof of Theorem \ref{main_408}:} Substitute $\eps =
n^{-1/14}$ in Theorem \ref{prop_408}. \hfill $\square$

\medskip {\emph {Remark:}} The exponents $3.33$ and $0.33$ in Theorem \ref{prop_408} are not optimal.
They may be replaced by constants arbitrarily close to $10/3$ and
$1/3$, respectively, at the expense of increasing $C$ and
decreasing $c$ in Theorem \ref{prop_408}, as is easily seen from
our proof. We conjecture that slightly better exponents may be
obtained by optimizing our argument; for example, the transition
from (\ref{improve1}) to (\ref{improve2}) seems inefficient, and
it also makes sense to try and play with the function $ M(U) =
\langle \nabla \log \pi_{U(E_0)} g(U x_0), U x_0 \rangle $ in
place of the definition (\ref{eq_1141}).

\section{Tying up loose ends}

Next we discuss the proof of Theorem \ref{main_630}. As in the
previous section, we rely heavily on results from \cite{CLT}. For
two random vectors $X$ and $Y$ attaining values in some measurable
space $\Omega$, we write
$$ d_{TV}(X, Y) = 2 \sup_{A \subseteq \Omega} \left|
\, Prob \{ X \in A \} \, - \, Prob \{ Y \in A \} \, \right|, $$
for the total-variation distance between $X$ and $Y$, where the
supremum runs over all measurable subsets $A \subseteq \Omega$.
The following lemma  is no more than an adaptation of
\cite[Proposition 5.7]{CLT}.

\begin{lemma}
Let $2 \leq \ell \leq n$ be integers and assume that $X$ is a
random vector in $\RR^n$ with an isotropic, log-concave density.
Suppose that $$ \ell \leq n^{\kappa}. $$

\smallskip Then, there exists a subset $\cE \subseteq G_{n,\ell}$
with $\sigma_{n,\ell} (\cE) > 1 - C e^{-\sqrt{n}}$ such that for
all $E \in \cE$ there exists a random vector $Y$ in $E$ for which
the following hold:
\begin{enumerate}
\item[(i)] $\displaystyle d_{TV}(Proj_E(X), Y) \leq C / \ell^{10} $.
\vspace{5pt} \item[(ii)]  $Y$ has a spherically-symmetric
distribution. That is, $$ Prob \{ Y \in A \} = Prob \{ Y \in U(A)
\} $$ for any measurable subset $A \subseteq E$ and any $U \in
SO(n)$ with $U(E) = E$.
 \vspace{3pt} \item[(iii)] $\displaystyle
 Prob \left \{ \big| \, |Y| - \sqrt{\ell} \, \big| \geq \eps
 \sqrt{\ell}
\right \} \leq C e^{-c \eps^2 \ell } $  for all $0 \leq \eps \leq
1$.
\end{enumerate}
Here, $c, C, \kappa > 0$ are universal constants.
 \label{prop_842}
\end{lemma}

\emph{Proof:} We may clearly assume that $n \geq \tilde{C}$. We
begin by verifying the requirements of Lemma \ref{lem_628}. Let
$\alpha = 3000, \eta = 10$ and $u = 1/3$. Our universal constant
$\kappa$ will be defined  by
$$ \kappa = \min \left
\{  1 / (4 \alpha + 2 \eta + 5.01)  ,  u / (4 \alpha + 2 \eta + 4)
\right \}.
$$
Recall that $\ell \leq n^{\kappa}$.  Let $f: \RR^n \rightarrow [0,
\infty)$ stand for the density of $X$, and denote $g = f *
\gamma_{n}[\ell^{- \alpha}]$. The requirements of Lemma
\ref{lem_628} thus hold true. By the conclusion of that lemma,
there exists $\cE \subseteq G_{n,\ell}$ with $\sigma_{n,
\ell}(\cE) > 1 - C \exp(-\sqrt{n})$ such that for $E \in \cE$,
\begin{equation}
\pi_E(g)(r \theta_1) \leq \left(1 + \frac{2}{\ell^{10}} \right)
\cdot \pi_E(g)(r \theta_2) \label{eq_619}
\end{equation}
for all $\theta_1, \theta_2 \in S(E), 0 \leq r \leq 10
\sqrt{\ell}$. Fix $E \in \cE$. We need to construct a random
vector $Y$ in $E$ that satisfies (i), (ii) and (iii). Consider
first the random vector $X^{\prime}$ in the subspace $E$ whose
density is $\pi_E(g) = \pi_E(f) * \gamma_E[\ell^{-\alpha}]$. The
function $\pi_E(f)$ is isotropic, log-concave and it is the
density of $Proj_E(X)$. According to \cite[Lemma 5.1]{CLT},
\begin{equation}
d_{TV} \left( Proj_E(X) ,  X^{\prime} \right) \leq
\frac{C}{\ell^{10}}. \label{eq_1031}
\end{equation}
The density $\pi_E(g)$ is $C^2$-smooth (it is a convolution with a
gaussian) and log-concave (by Pr\'ekopa-Leindler). Additionally,
\begin{equation}
 \ell \leq \EE |X^{\prime}|^2 = \ell + \ell^{1- \alpha} \leq 2
\ell. \label{eq_1032}
\end{equation}
 We may thus apply Lemma \ref{lem_853}, based  on (\ref{eq_619}), for
 $\delta = 2 / \ell^{10}$. According to the
conclusion of Lemma \ref{lem_853},
\begin{equation} Prob \left \{ \left| \frac{|X^{\prime}|}
{\sqrt{\ell + \ell^{1-\alpha}}} - 1  \right| > \eps \right \} \leq
C e^{-c \eps^2 \ell } \ \ \ \text{for all} \ 0 \leq \eps \leq 1.
\label{eq_1029}
\end{equation}
 Since $\alpha > 1/2$ then
$\sqrt{\ell + \ell^{1-\alpha}}$ is sufficiently close to
$\sqrt{\ell}$, and from (\ref{eq_1029}) we obtain
\begin{equation} Prob \left \{ \left| \frac{|X^{\prime}|}
{\sqrt{\ell}} - 1  \right| > \eps \right \} \leq C^{\prime}
e^{-c^{\prime} \eps^2 \ell } \ \ \ \text{for all} \ 0 \leq \eps
\leq 1. \label{eq_1033}
\end{equation}
Define $Y$ to be the random vector in the subspace $E$ whose
density is
\begin{equation}
 \tilde{g}(x) = \int_{S(E)} \pi_E(g)(|x| \theta)
d\sigma_{E}(\theta). \label{eq_1159}
\end{equation}
 Then $\tilde{g}$ is the spherical average
of $g$, hence $Y$ is spherically-symmetric in $E$ and (ii) holds.
Additionally, since $|X^{\prime}|$ and $|Y|$ have the same
distribution, then (iii) holds in view of (\ref{eq_1033}). All
that remains is to prove (i). According to (\ref{eq_619}) and
(\ref{eq_1159}), for any $x \in E$ with $|x| \leq 10 \sqrt{\ell}$,
\begin{equation}
 |\tilde{g}(x) -
 \pi_E(g)(x)| \leq \frac{C}{\ell^{10}} \cdot \tilde{g}(x).
\label{eq_1122}
\end{equation}
From (\ref{eq_1033}), (iii) and (\ref{eq_1122}),
\begin{eqnarray}
\nonumber \lefteqn{d_{TV}(Y, X^{\prime}) = \int_{\RR^n}
|\tilde{g}(x) -  \pi_E(g)(x) | dx } \\ & \leq & \nonumber 4 Prob
\{ |X^{\prime}| \geq 2 \sqrt{\ell} \} + \int_{|x| \leq 2
\sqrt{\ell}} |\pi_E(g)(x) - \tilde{g}(x)| dx \\ & \leq & C e^{-c
\ell} + \frac{C^{\prime}}{\ell^{10}} \int_{|x| \leq 2 \sqrt{\ell}}
\tilde{g}(x) dx \leq \frac{\tilde{C} }{\ell^{10}}. \label{eq_1127}
\end{eqnarray}
Now (i) follows from (\ref{eq_1127}) and (\ref{eq_1031}). This
completes the proof. \hfill $\square$

\medskip \emph {Proof of Theorem \ref{main_630}:} We will actually prove
the theorem under the weaker assumption that the density of $X$ is
log-concave. The deduction of Theorem \ref{main_630} from Lemma
\ref{prop_842} is very similar to the argument in the proof of
Theorem 5.9 in \cite{CLT}. We supply a few details. We may assume
that $n \geq \tilde{C}$, since otherwise there is no $\ell \geq 1$
with $\ell \leq c n^{\kappa}$. It is sufficient to prove the
theorem for $\ell = \lfloor n^{\kappa / 2} \rfloor + 2$ where
$\kappa$ is the constant from Lemma \ref{prop_842}. Let $\cE
\subseteq G_{n, \ell^2}$ be the subset from the conclusion of
Lemma \ref{prop_842} for $\ell^2$. Then
\begin{equation}
\sigma_{n, \ell^2}(\cE) \geq 1 - C \exp(-\sqrt{n}).
\label{eq_1147}
\end{equation}
Fix $E \in \cE$, and let $Y$ be the random vector in $E$ whose
existence is guaranteed by Lemma \ref{prop_842}. Let $F \subset E$
be any $\ell$-dimensional subspace. We may apply \cite[Lemma
5.8]{CLT}, based on properties (ii) and (iii) of $Y$ from Lemma
\ref{prop_842}. We conclude that
\begin{equation}
 d_{TV} \left( Proj_F(Y)  , Z_F \right) \leq C
\sqrt{\frac{\ell}{\ell^2}} \leq \frac{C^{\prime}}{\sqrt \ell},
\label{eq_1138_}
\end{equation}
 where $Z_F$ is a standard gaussian
random vector in the subspace $F$. Recall that $d_{TV}(Y,
Proj_E(X)) < C / \ell^{10}$, by property (i) from Lemma
\ref{prop_842}. Therefore, from (\ref{eq_1138_}),
\begin{equation}
d_{TV} \left( Proj_F(X)  , Z_F \right) \leq
\frac{C^{\prime}}{\sqrt{\ell}} + \frac{C}{\ell^{10}} \leq
\frac{\bar{C}}{n^{\kappa / 4}}. \label{eq_1138__}
\end{equation}
Recall that $E \in \cE$ and $F \subset E$ were arbitrary. Denote
$$ \cF =\left\{ F \in G_{n,\ell } ; \exists E \in \cE, \ F
\subset E \right \}. $$ We have proved that (\ref{eq_1138__})
holds for all $F \in \cF$. From (\ref{eq_1147}) we deduce that
$\sigma_{n, \ell}(\cF) \geq \sigma_{n, \ell^2}(\cE) \geq 1 -
\exp( - c \sqrt{n} )$. The theorem is thus proven. \hfill
$\square$

\medskip The remaining statements that were announced in Section
\ref{intro}  above follow, in a rather straightforward manner,
from the theorems obtained so far in this note and from results
found in the literature. The argument that leads from Theorem
\ref{main_630} to Corollary \ref{cor_847} is elementary and
well-known (in the context of Dvoretzky's theorem, we learned it
from G. Schechtman). It is based on the observation that any
$(2n+1)$-dimensional ellipsoid $\cE \subset \RR^{2n + 1}$ has a
projection onto some $n$-dimensional subspace $F \subset \RR^{2n +
1}$ such that $Proj_F(\cE)$ is a multiple of the standard
Euclidean ball in the subspace $F$. We omit the standard linear
algebra details.

\medskip Regarding the proof of Theorem \ref{thm_sodin}: We may assume
that $n \geq C$. Note that the desired conclusion (i) is
equivalent to the case $\ell = 1$ in Theorem \ref{main_630}, since
the total-variation distance between two random variables equals
the $L^1$ distance between their densities. In order to prove
(ii), we use Sodin's results \cite{sodin}. We may apply
\cite[Theorem 2]{sodin} with $\alpha = 0.33$ and $\beta = 3.33$,
in view of Theorem \ref{prop_408}. According to the conclusion of
\cite[Theorem 2]{sodin}, for any $t \in \RR$ with $|t| < c
n^{1/24}$,
\begin{equation}
 \left| \int_{S^{n-1}} [f_{\theta}(t) / \gamma(t)]
d\sigma_{n-1}(\theta) - 1 \right| \leq \frac{C}{n^{1/24}}.
\label{eq_950}
\end{equation}
 Next, we would like to use a ready-made concentration of measure phenomenon
 argument, such as \cite[Theorem 5]{sodin}. However, the results in the
 relevant section in \cite{sodin} are proven under the additional
 assumption that $X$ is a symmetric random vector
(i.e., $X$ and $-X$ have the same distribution). Sodin's argument
 formally relies on the fact that in the symmetric case, $f_{\theta}(t)$ is
non-increasing for $t > 0$, as an even, log-concave density. In
the case where $X$ is symmetric, we may directly apply
\cite[Theorem
 5]{sodin} for $\eps = C n^{-1/24}$ and $T = c n^{1/24}$, because of (\ref{eq_950}).
 We deduce that there exists $\Theta \subseteq S^{n-1}$  with
 $\sigma_{n-1}(\Theta) \geq 1 - C \exp (-\sqrt{n})$ such that for
 any $\theta \in \Theta$,
 $$ |f_{\theta}(t) / \gamma(t) - 1| \leq
 \frac{C^{\prime}}{n^{1/24}} \ \ \ \text{when} \ \ |t| < c n^{1 /
 24}. $$
This completes the proof of Theorem \ref{thm_sodin} for the case
where the density of $X$ is an even function.

\smallskip We  claim that Theorem \ref{thm_sodin}(ii) is true as stated,
without the additional assumption that the random vector $X$ is
symmetric; it is possible to modify Sodin's argument
(specifically, Proposition 13 in \cite{sodin}) for the general,
log-concave case. We will not carry out the details here, and they
are left to the reader as a (rather interesting) exercise. An
alternative route to establish Theorem \ref{thm_sodin}(ii) from
(\ref{eq_950}), in the general, non-even case, may be very roughly
summarized as follows: Observe that after a convolution with a
small gaussian, estimates such as (\ref{eq_619}) directly lead us
to the desired result. Then, show that the convolution of a
log-concave function with a small gaussian has only a minor effect
in the moderate-deviation scale.

\medskip \emph{Remark:} It is also possible to improve the
quantitative bound for $\eps_n$ from \cite[Theorem 1.2]{CLT}. The
most straightforward adaptation of the proof of Theorem 1.2 in
\cite{CLT}, using the new Theorem \ref{prop_408}, leads to the
estimate $\eps_n \leq C / n^{\kappa}$ for some universal constants
$C, \kappa > 0$.

\bigskip
{\small  \noindent Department of Mathematics, Princeton
university, Princeton, NJ 08544, USA \\ {\it e-mail address:}
\verb"bklartag@princeton.edu"

\end{document}